\documentclass{amsart}
\usepackage{amssymb}
\usepackage{latexsym}
\newtheorem{thm}{Theorem}

\newtheorem{asp}{Assumption}
\newtheorem{lem}{Lemma}
\newtheorem{cor}{Corollary}
\newtheorem{prop}{Proposition}
\newtheorem{dfn}{Definition}
\begin{document}
\title{Volume minimization for Lagrangian submanifolds in complex manifolds with negative first Chern class}
\author{Edward Goldstein}
\maketitle
\renewcommand{\abstractname}{Abstract}
\begin{abstract}
It is classically known that closed geodesics on a compact Riemann surface with a metric of negative curvature strictly minimize length in their free homotopy class. We'd like to generalize this to Lagrangian submanifolds in K\"ahler manifolds of negative Ricci curvature. The only known result in this direction is a theorem on Y.I. Lee for certain Lagrangian submanifolds in a product of two Riemann surfaces of constant negative curvature. We develop an approach to study this problem in higher dimensions. Along the way we prove some weak results (volume-minimization outside of a divisor) and give a counterexample to global volume-minimization for an immersed minimal Lagrangian submanifold.
\end{abstract}
\section{Local theory}
In this section we'll present local volume-minimizing properties for minimal Lagrangian submanifolds in K\"ahler manifolds with negative Ricci curvature. Strictly speaking the rest of the paper does not logically depend on this section and this section serves as a motivation for the global problem.\\
Let $(M^{2n},\omega)$ be a K\"ahler $n$-fold and let $K(M)$ be the canonical bundle of $M$, i.e. $K(M)=\Lambda^{(n,0)}T^{\ast}M$- the bundle of $(n,0)$-forms on $M$. So $c_1(K(M))=-c_1(M)$. The Hermitian metric on $K(M)$ endowed from $M$ gives rise to a Hermitian connection on $K(M)$. The curvature form of $K(M)$ is called the Ricci form of $M$ and denoted by $Ric$. It is related to the Ricci tensor of $M$ (which is a symmetric 2-tensor on $M$) by the following identity:
\begin{equation}
\label{Ricci curvature vs. tensor}
-iRic(X,Y)=Ricci(JX,Y)
\end{equation}
for any tangent vectors $X$ and $Y$ (here $J$ is the complex structure on $M$).\\
Let $L$ be an oriented Lagrangian submanifold of $M$. For any point $l \in L$ there is a unique element $\kappa_l$ of $K(M)$ over $l$ which restricts to the volume form on $L$ (in the induced metric from $M$). Various $k_l$ give rise to a section 
\begin{equation}
\label{can}
\kappa_L:L \rightarrow K(M)
\end{equation}
The section $\kappa_L$ has constant length $\sqrt{2}^n$ over $L$ and it
defines a connection 1-form $\xi$ on $L$ for the connection on $K(M)$ over $L$ by the condition $\xi \otimes \kappa_L = \nabla \kappa_L$.\\
Let $H$ be the mean curvature vector field of $L$. So $H$ is a section of the
normal bundle of $L$ in $M$ and we have a corresponding 1-form $\sigma=
i_H\omega$ on $L$. The following fact was first observed by Y.G. Oh in \cite{Oh1}, see also \cite{CG}:
\begin{lem} 
\label{mc}
$i\sigma= \xi$
\end{lem}
Thus if $L$ is {\it minimal} (i.e. $H=0$) iff $\kappa_L$ is parallel over $L$.
From now on we work with $L$ which is minimal, embedded and compact. We''ll also assume that $L$ is a real-analytic submanifold of $M$ - this is always true if the K\"ahler metric is real-analytic. Since $L$ is a totally-real real-analytic submanifold of $M$ we easily see that we can extend $\kappa_L$ to be a holomorphic section of $K(M)$ in some neighbourhood $U$ of $L$ in $M$. Consider the function $f=ln|\kappa_L|^2$ on $U$. $L$ is a critical submanifold of $f$ on $U$. Moreover $K(M)$ is a positive line bundle. It follows that at a point $l \in L$ the Hessian $Hess(f)$ of $f$ on $M$ has at least $n$ negative eigenvalues - see \cite{GH}, p. 158. Thus the Hessian of $f$ is negative-definite on the normal bundle to $M$ on $L$. Hence $f$ and $|\kappa_L|$ attains a strict local maximum along $L$. Moreover $\kappa_L$ restricts to the volume form on $L$. Thus $\kappa_L$ calibrates $L$ (see \cite{HL}).
\begin{lem}
$L$ is calibrated by $Re\kappa_L$ in some neighbourhood $U'$ of $L$ in $M$ and moreover $L$ is {\bf strictly} homologically volume-minimizing in $U'$.
\end{lem}
We would like to globalize this statement from a neighbourhood of $L$. If $M$ is a compact Riemann surface then the following result is classical:
\begin{lem}
On a compact Riemann surface with a metric of negative curvature every free homotopy class of loops contains a unique geodesic and this geodesic strictly minimizes length in this class.
\end{lem}
The only progress attained so far in higher dimensions is the result of Y.I. Lee [Lee] in case $M$ is a product of two Riemann surfaces with constant negative curvature. For instance  Y.I. Lee has shown that a product of two geodesics is {\bf homotopically} area minimizing. We'll now describe a program of studying this issue in higher dimensions.
\section{Algebraic minimal Lagrangian submanifolds}
We continue to consider a K\"ahler $n$-fold $M$ with negative Ricci curvature. As a first step we'll study some ``nice'' minimal Lagrangian submanifolds, which we'll call {\it algebraic}:
\begin{dfn}
Let $L$ be a minimal Lagrangian submanifold of $M$. We say that $L$ is {\it algebraic} if there is a holomorphic section $\varphi$ of $K(M)$ over $M$ s.t. $|\varphi|$ has a strict global maximum along $L$.
\end{dfn}
\begin{lem}
Algebraic minimal Lagrangian submanifolds strictly minimize volume in their homology class.
\end{lem}
{\bf Proof:} Note that since $|\varphi|$ has a local maximum along $L$,  $\varphi$ is parallel along $L$. Since $L$ is minimal it follows from Lemma \ref{mc} that we can multiply $\varphi$ by a unit complex number $e^{i\theta}$ so that $e^{i\theta}\varphi$ will restrict to the volume form on $L$. Thus $L$ is calibrated by $e^{i\theta}\varphi$. Q.E.D.\\
Next we'll demonstrate some examples of algebraic minimal Lagrangian submanifolds. Let $M$ be a complete intersection in $\mathbb{C}P^n$ given by polynomials $f_1,\ldots,f_m$ with real coefficients. We want $\Sigma deg(f_i)-(n+1)=2k>0$ for an integer $k$. Then $K(M) \simeq (\gamma^{\ast})^{\otimes2k}$ where $\gamma^{\ast}$ is the hyperplane bundle of $\mathbb{C}P^n$. Let $L=N \cap \mathbb{R}P^n$. $L$ is a minimal Lagrangian submanifold.\\
Consider the standard Hermitian metic on $(\gamma^{\ast})^{\otimes 2k}$. We have a section \[\varphi=(\Sigma z_i^2)^k\]
of $(\gamma^{\ast})^{\otimes2k}$ and $|\varphi|$ attains strict maximum along $\mathbb{R}P^n$. By Yau's resolution of the Calabi conjecture [Yau] we can put a K\"ahler metric $\omega_Y$ on $M$ s.t. the induced metric on $K(M)$ will be proportional by a constant to the standard metric on $(\gamma^{\ast})^{\otimes2k}$. Thus we'll get:
\begin{prop}
$L$ is an algebraic minimal Lagrangian submanifold of $(M,\omega_Y)$.
\end{prop}
\section{A global strategy} 
We''ll now present our strategy for exploring global volume-minimizing properties for Lagrangian submanifolds in complex manifolds with positive canonical bundle. We'll explain this approach in the simplest possible case - a fixed point set of an antiholomorphic involution. This approach can be modified to a more general setting - see the remark in the end of Section \ref{berg}.\\
So let $M$ be a complex $n$-fold with an antiholomorphic involution $\sigma$. Let $g$ be a Hermitian metric on $M$ which is $\sigma$-invariant. Note we don't assume that $g$ is K\"ahler, just Hermitian. 
\begin{asp}
The metric $g$ induces a Hermitian metric $h$ on the canonical bundle $K(M)$ of the holomorphic $(n,0)$-forms of $M$. Let $\omega$ be the curvature form of this metric. We assume that $\omega$ defines a K\"ahler form on $M$. 
\end{asp}
If $g$ were a K\"ahler metric this would be equivalent to $g$ having negative Ricci curvature. Let $L$ be the fixed pont set of $\sigma$. We'd like to understand when $L$ minimizes its $g$-volume in its free homotopy class. First we'll describe the convergence of Bergmann metrics following \cite{Tian} and \cite{Ruan}.
\subsection{Bergman metrics}
\label{berg}
Let $H^0(K(M)^{\otimes m},M)$ be the space of holomorphic sections of the $m^{th}$ power $K(M)^{\otimes m}$. Clearly $\sigma$ acts upon $H^0(K(M)^{\otimes m},M)$ by an antiholomorphic involution (the $\sigma$-action is by pulling back a section and conjugating it). Let $H_m^L$ be the fixed point set of $\sigma$ on $H^0(K(M)^{\otimes m},M)$. Then $H_m^L$ is a Lagrangian subspace of $H^0(K(M)^{\otimes m},M)$. Pick an orthonormal basis $\kappa_i$ for $H_m^L$ - here the volume form on $M$ is given by $\frac{\omega^n}{n!}$ and this gives a natural inner product on $H^0(K(M)^{\otimes m},M)$. Clearly $\kappa_i$ it will give a unitary basis for $H^0(K(M)^{\otimes m},M)$. It gives rise to a Kodaira embedding: \[\phi\sb {m}: M \mapsto \mathbb{C}P\sp {N\sb {m}}\]
Moreover $\phi\sb {m}$ will be intertwining between $\sigma$ and the conjugation on $\mathbb{C}P\sp {N\sb {m}}$. Thus \[L= \phi\sb {m}^{-1}(\mathbb{R}P\sp {N\sb {m}})\]
Let $\omega_{FS}$ be the Fubini-Study metric on $\mathbb{C}P\sp {N\sb {m}}$ and let \[\omega_m=(\phi\sb {m})^{\ast}\omega_{FS}/m\]
be the corresponding Bergmann metric on $M$. Then by \cite{Tian} and \cite{Ruan}, $\omega_m \rightarrow \omega$ on $M$. Here the convergence is $C^{\infty}$.\\
The Kodaira embedding $\phi\sb {m}$ also gives an isomorphism between $K(M)^{\otimes m}$ and the hyperplane bundle $\gamma^{\ast}$ of $\mathbb{C}P\sp {N\sb {m}}$. How $\gamma^{\ast}$ has a natural Hermitian metric and so it induces a Hermitian metric on $K(M)^{\otimes m}$ and hence on $K(M)$. Let $h_m$ be that metric on $K(M)$. Then the curvature form of $h_m$ is $\omega_m$. We'd like to conformally scale $g$ so that the induced metric on $K(M)$ will equal to $h_m$ up to a constant multiple. If the scaling is $g_m=e^{f_m}g$ then $f_m$ solves:
\[ n i\partial \bar{\partial}f_m=\omega_m-\omega \]
Here $n$ is the dimension of $M$. If we choose $f_m$ to intergate to $0$ over $M$ then we get that $f_m \rightarrow 0$ in the $C^{\infty}$-topology. Hence 
\begin{equation}
\label{metcon}
g_m \stackrel{C^{\infty}}{\rightarrow} g 
\end{equation}
Consider now a holomorphic section $\Sigma z_i^2$ of $(\gamma^{\ast})^{\otimes 2}$ over $\mathbb{C}P\sp {N\sb {m}}$. We pull it back to a holomorphic section
$\varphi_m$ of $K(M)^{\otimes 2m}$ over $M$. In fact
\[\varphi_m=\Sigma \kappa_i^{\otimes 2}\]
Here $\kappa_i$ form an orthonormal basis of $H_m^L$. The upshot of our construction is:
\begin{prop}
The section $\varphi_m$ of $K(M)^{\otimes 2m}$ is $\sigma$-invariant and its length with respect to $g_m$ attains a maximum along $L$.
\end{prop}
{\bf Remark:} We have a canonical construction of the sections $\varphi_m$ for $L$ being the fixed point set of an antiholomorphic involution. There is a more general (non-canonical) construction with similar (somewhat weaker) properties using Donaldson's asymptotically holomorphic sections, see \cite{AGM}.
\subsection{Homotopy volume minimization outside of a divisor}
We now would like to use the idea of calibration to study the $g_m$-volume minimizing properties of $L$. We have the following result:
\begin{thm}
Suppose we have a 1-parameter family $L_t$ of homotopies of $L$ with $L_0=L$ and moreover all $L_t$ avoid the zero set $D_m=\varphi_m^{-1}(0)$. Then $vol(L_t) \geq vol(L)$. Here the volume is with respect to $g_m$.
\end{thm}
{\bf Proof:} First will prove this for $L$ orientable. For non-orientable $L$ we look on its orientable double cover and copy the proof bellow.\\
Since $\varphi_m$ is $\sigma$-invariant we have choose its $2m^{th}$ root $\xi_m$ along $L$ such that $\xi_m$ restricts to a positive real $n$-form on $L$. We extend $\xi_m$ to be a $2m^{th}$ root of $\varphi_m$ over $L_t$ along the homotopy. Then $\xi_m$ has constant length $C$ (with respect to $g_m$) along $L$ and the length of $\xi_m$ is $\leq C$ on $L_t$. We obviously have \[\int_{L} \xi_m= \int_{L_t} \xi_m\] and hence $vol(L_t) \geq vol(L)$.   Q.E.D.\\
The last statement about the conservation of the integral of $\xi_m$ can be better understood using the ramified $2m$-cover of $M$ that $\varphi_m$ defines. The setup is as follows:\\ 
Let $K(M)$ be the total space of the canonical bundle of $M^{2n}$ and $\pi : K(M) \rightarrow M$ be the projection. There is a canonical $(n,0)$-form $\rho$ on $K(M)$ defined by $\rho(a)(v_1,\ldots,v_n)=a(\pi_{\ast}(v_1),\ldots,\pi_{\ast}(v_n))$. Here  $a \in K(M)$ and $v_1,\ldots,v_n$ are tangent vectors to $K(M)$ at $a$. The form $\eta = d\rho$ is a holomorphic volume form on $K(M)$. If $z_1, \ldots,z_n$ are local coordinates on $M$ then $dz_1\wedge \ldots \wedge dz_n$ is a local section of $K(M)$ over $M$ which defines a coordinate function $y$ on $K(M)$. The collections of holomorphic functions $(z_1,\ldots,z_n,y)$ are coordinates on $K(M)$ and  
\begin{equation}
\label{rho-eta}
\rho=ydz_1 \wedge \ldots \wedge dz_n ~ , ~ \eta=dy \wedge dz_1\wedge \ldots \wedge dz_n
\end{equation}
Now we have a section $\varphi_m$ of $K(M)^{\otimes 2m}$ and it defines a $2m$-fold ramified cover $M_m$ of $M$ as follows:
\[M_m=(\alpha \in K(M)| \alpha^{\otimes 2m}=\varphi_m(\pi(\alpha)))\]
We can in fact slightly perturb $\varphi_m$ to keep it holomorphic, $\sigma$-invariant and make it transversal to $0$. Than way $M_m$ is a complex $n$-fold and the projection $\pi:M_m \rightarrow M$ is a ramified cover. Moreover the $(n,0)$-form $\rho$ on $K(M)$ restricts to a holomorphic form on $M_m$.\\ 
Now the $2m^{th}$ root $\xi_m$ of $\varphi_m$ along $L$ can be viewed to a lift $L^K$ of $L$ to $M_m$. If we have a homotopy $L_t$ of $L$ in $M$ avoiding the zero set $N_m=\varphi_m^{-1}(0)$ then $L_t$ can be lifted to a homotopy $L_t^K$ on $M_m$. Obviously we have \[ \int_{L_t^K}\rho= \int_{L^K}\rho\]
Hence the volume comparison property. We have a stronger result:
\begin{cor}
Let $L'$ be a submanifold of $M_m$ representing the same homology class as $L^K$ in $M_m$. Then $vol(\pi(L')) \geq vol(L)$. Here again the volume is $g_m$-volume.
\end{cor}
We'd like eventually to understand when $L$ minimizes $g$-volume in its free homotopy class. Note that if $L$ minimizes $g$-volume then it strictly minimizes $g$-volume. Indeed suppose that some other $L'$ homotopic to $L$ has the same $g$-volume. Pick a point $a$ in $L$ and not in $L'$. Consider a positive $\sigma$-invariant function $f$ supported in a neighbourhood of $a$ and consider a Hermitian metric $g_f=e^fg$. Then the $g_f$ volume of $L$ is greater then the $g_f$-volume of $L'$. We can also choose $f$ small enough so that the Hermitian metric that $g_f$ induces on $K(M)$ still has positive curvature and this will lead us to a contradiction. In particular we'll deduce that if $\alpha$ is an automorphism of $M$ then $L$ and $\alpha(L)$ are not homotopic if they don't coincide.
\section{A null-homotopic immersed minimal Lagrangian sphere}
In this section we'll give an example of a compact K\"ahler-Einstein surface $M$ with negative Ricci curvature and an immersed minimal Lagrangian sphere $i: S^2 \mapsto M$ which is null-homotopic. \\
The example is: $M=(\Sigma z_i^5=0 \subset \mathbb{C}P^3)$ endowed with a K\"ahler-Einstein metric. Let $L$ be a fixed point set of the conjugation on $M$. Then $L$ is clearly diffeomorphic to $\mathbb{R}P^2$. Let $S^2$ be its orientable double cover. Thus we have an immersion $i:S^2 \mapsto M$. Clearly $i(S^2)$ is a $2$-torion in $H_2(M,\mathbb{Z})$ and since $M$ is a simply connected 4-manifold, $H_2(M,\mathbb{Z})$ has no torsion. Hence $i(S^2)$ is trivial in $H_2(M,\mathbb{Z})$ and by the Hurewicz theorem $i(S^2)$ is null-homotopic. \\
From the results of J. Lees \cite{Lees} we in fact can also deduce that this Lagrangian immersion of $S^2$ can be homotoped through Lagrangian immersions to a Whitney sphere contained in a Darboux chart of some point in $M$. We omit the details.\\
On the other hand $L$ itself is non-trivial in the second homology of $M$ with $\mathbb{Z}/2$-coefficients. Moreover the results of \cite{Gol} show the following: view $M$ with the metric induced from the Fubini-Study metric on $\mathbb{C}P^3$. Let $L'$ be a Lagrangian submanifold of $M$ in the same second $\mathbb{Z}/2$-homology class as $L$. Then \[5 \cdot Area(L') \geq Area(L) \] 
We hope that global volume-minimizing properties are true for embedded minimal Lagrangian submanifolds.

Institute for Advanced Study, Princeton NJ\\
egold@ias.edu
\end{document}